\documentclass[preprint,noinfoline]{imsart}

\usepackage{amsmath,amsthm,amsfonts,amssymb}
\usepackage[authoryear,longnamesfirst]{natbib}
\usepackage{enumitem}
\usepackage{color}
\usepackage[pdftitle={The Contact Process in a Dynamic Random Environment},pdfauthor={Daniel Remenik}]{hyperref}
\usepackage{hypernat}

\definecolor{darkblue}{rgb}{0.03,0.03,0.23}
\hypersetup{colorlinks=true,urlcolor=darkblue,citecolor=darkblue,linkcolor=darkblue}

\startlocaldefs

\newtheorem{teo}{Theorem}
\newtheorem{teosec}{Theorem}[section]
\newtheorem{lem}[teosec]{Lemma}
\newtheorem{prop}[teosec]{Proposition}
\newtheorem{cor}[teosec]{Corollary}
\theoremstyle{definition}
\newtheorem{rem}[teosec]{Remark}

\newcommand{\pp}{\mathbb{P}}
\newcommand{\ee}{\mathbb{E}}
\newcommand{\rr}{\mathbb{R}}
\newcommand{\nn}{\mathbb{N}}
\newcommand{\zz}{\mathbb{Z}}

\newcommand{\cf}{\mathcal{F}}
\newcommand{\cx}{\mathcal{X}}

\newcommand{\mr}{\mu_\rho}
\newcommand{\nA}{{\nu_A}}
\newcommand{\nC}{{\nu_C}}

\newcommand{\notempty}{\neq\emptyset}

\newcommand{\onu}{\ensuremath{\overline{\nu}}}
\newcommand{\unu}{\ensuremath{\underline{\nu}}}
\newcommand{\dual}[2]{\hat{#1}^{#2}}
\newcommand{\uno}[1]{\mathbf{1}_{#1}}

\newcommand{\BC}{\hyperlink{BC}{{\rm(BC)}}}

\let\oldmarginpar\marginpar
\renewcommand\marginpar[1]{\oldmarginpar{\raggedright\texttt{#1}}}

\numberwithin{equation}{section}

\endlocaldefs

\begin{document}

\begin{frontmatter}

  \title{The Contact Process in a Dynamic Random Environment} \runtitle{The Contact
    Process in a Dynamic Random Environment}

  \author{\fnms{Daniel} \snm{Remenik}\ead[label=e1]{dir4@cornell.edu}}

  \affiliation{Cornell University} \runauthor{D. Remenik}

  \address{Center for Applied Mathematics\\
    Cornell University\\
    657 Rhodes Hall,\\
    Ithaca, New York 14853\\
    USA\\
    \printead{e1}}

\begin{abstract}
  We study a contact process running in a random environment in $\zz^d$ where sites flip,
  independently of each other, between blocking and non-blocking states, and the contact
  process is restricted to live in the space given by non-blocked sites. We give a partial
  description of the phase diagram of the process, showing in particular that, depending
  on the flip rates of the environment, survival of the contact process may or may not be
  possible for large values of the birth rate. We prove block conditions for the process
  that parallel the ones for the ordinary contact process and use these to conclude that
  the critical process dies out and that the complete convergence theorem holds in the
  supercritical case.
\end{abstract}



\end{frontmatter}

\section{Introduction}\label{sec:intr}

We consider the following version of a contact process running in a dynamic random
environment in $\zz^d$. The state of the process is represented by some
$\eta\in\cx={\{-1,0,1\}}^{\zz^d}$, where sites in state 0 are regarded as vacant, sites in
state 1 as occupied, and sites in state $-1$ as blocked (that is, no births of 1's are
allowed on that site). The process $\eta_t$ is defined by the following transition rates:
\begin{center}
  \begin{tabular}{ccclc}
    $0$ & $\longrightarrow$ & $1$\hspace{.3in} & at rate & $\beta f_1$ \\
    $1$ & $\longrightarrow$ & $0$\hspace{.3in} & at rate & 1 \\
    $0,1$ & $\longrightarrow$ & $-1$\hspace{.3in} & at rate & $\alpha$ \\
    $-1$ & $\longrightarrow$ & $0$\hspace{.3in} & at rate & $\alpha\delta$ \\
  \end{tabular}
\end{center}
where $f_1$ is the fraction of occupied neighbors at $L^1$ distance 1.

In words, the $-1$'s define a random environment in which each site becomes blocked at
rate $\alpha$ and flips back to being unblocked at rate $\alpha\delta$, while the 1's
behave like a nearest neighbor contact process with birth rate $\beta$ in the space left
unblocked by the environment. Observe that when an occupied site becomes blocked, the
particle is killed. This version is simpler than the alternative in which only 0's can
turn to $-1$'s (mainly because our process satisfies a self-duality relation, see
Proposition \ref{prop:dual}). However, we feel that our choice is natural: if a site
becomes uninhabitable, the particles living there will soon die.

Ever since it was introduced in \citet*{harris1}, the contact process has been object of
intensive study, and many extensions and modifications of the process have been
considered. In particular, the literature includes several different versions of contact
processes in random environments. One class of these processes corresponds to contact
processes where the birth and death rates are not homogeneous in space, and they are
chosen according to some probability distribution, independently across sites, and remain
fixed in time (see, for example, \citet{bramDurrSchonmann}, \citet{ligg3}, \citet{andjel},
and \citet{klein}). The main question for this class of processes is to determine
conditions on the parameters that guarantee or preclude survival.

A different class of models, which are somehow closer to the process we consider, have two
species with different parameters or ranges, but one of them behaves independently of the
other while the second is restricted to live in the space left by the first. These
processes were studied in \citet{durrettSwindle}, \citet{durrettMoller}, and
\citet{durrettSchinazi}. The results in these papers (mainly bounds on critical parameters
for coexistence and complete convergence theorems) are asymptotic, in the sense that they
are proved when the range of one or both types is sufficiently large.

The process we consider differs from both of the types of examples mentioned above: the
random environment is dynamic, and it behaves independently across sites. An example of a
spin system running in this type of environment was studied in \citet{luo}, and
corresponds to the Richardson model which would result from ignoring transitions from 1 to
0 in our process. Another example was studied recently in \citet{broman}, where the author
considers a process in which the environment changes the death rate of the contact process
instead of blocking sites. The dynamics of the process
$\Psi^{\gamma,p,A}_{\delta_0,\delta_1}$ introduced there are the same as those of our
process if $\delta_1=\infty$. The author considers this case as a tool in the study of his
process, but the results of the paper focus on the case $\delta_1<\infty$. We will use one
of his results to give a bound on a part of the phase diagram of our process in Theorem
\ref{thm:parameters}.

As mentioned above, the $-1$'s evolve independently of the 1's. They follow an
``independent flip process'' whose equilibrium is given by the product measure
\[\mr(\{\eta:\,\eta(x)=-1\})=1-\mr(\{\eta:\,\eta(x)\neq-1\})=\rho=\frac{1}{1+\delta}\qquad\forall
x\in\zz^d.\] This process is reversible, and its reversible measure is given by $\mr$.

In Section \ref{subsec:constMon} we will construct our process using the so-called
graphical representation. A direct consequence of the construction will be that $\eta_t$
satisfies some monotonicity properties analogous to those of the contact process. (Here
and in the rest of the paper, when we refer to the contact process we mean the
``ordinary'' nearest-neighbor contact process in $\zz^d$).  We consider the following
partial order on configurations:
\begin{equation}\label{eq:order}
  \eta^1\leq\eta^2\quad\Leftrightarrow\quad\eta^1(x)\leq\eta^2(x)\quad\forall
  x\in\zz^d.
\end{equation}
With this order, our process has the following property: given two initial states
$\eta^1_0\leq\eta^2_0$, it is possible to couple two copies of the process $\eta^1_t$ and
$\eta^2_t$ with these initial conditions in such a way that $\eta^1_t\leq\eta^2_t$ for all
$t\geq0$. We will refer to this property as \emph{attractiveness} by analogy with the case
of spin systems (this property is sometimes termed \emph{monotonicity}, see Sections II.2
and III.2 of \citet{ligg1} for a discussion of general monotone processes and of
attractive spin systems, respectively).

For $A\subseteq\zz^d$ we define the following probability measure $\nA$ on $\cx$: $-1$'s
are chosen first according to their equilibrium measure $\mr$ and then 1's are placed at
every site in $A$ that is not blocked by a $-1$. These measures are the initial conditions
for $\eta_t$ that are suitable for duality.

Let $\unu=\nu_\emptyset$, which corresponds to having the $-1$'s at equilibrium and no
1's. Let also $\onu$ be the limit distribution of the process when starting at the
configuration having all sites at state 1, which is obviously the largest configuration in
the partial order \eqref{eq:order}. We will show in Proposition \ref{prop:mon} that this
limit is well defined and it is stationary, and that $\unu$ and $\onu$ are, respectively,
the \emph{lower} and \emph{upper invariant measure} of the process (that is, the smallest
and largest stationary distribution of the process).

We will say that the process \emph{survives} if there is an invariant measure $\nu$ such
that
\[\nu\!\left(\left\{\eta\!:\,\eta(x)=1\text{ for some }x\in\zz^d\right\}\right)>0,\] or,
equivalently, if $\unu\neq\onu$ (we remark that, as a consequence of Theorem
\ref{thm:compConv}, every invariant measure for the process is translation invariant, so
the above probability is actually 1 whenever it is positive). Otherwise, we will say that
the process \emph{dies out}. We will see in Section \ref{sec:block} that this definition
of survival is equivalent to the following condition: the process started with a single 1
at the origin and everything else at $-1$ contains 1's at all times with positive
probability.

A second monotonicity property that will follow from the construction of $\eta_t$ is
monotonicity with respect to the parameters $\beta$ and $\delta$:
\begin{enumerate}[label=(\roman{*})]
\item If $\alpha$ and $\delta$ are fixed, and for some $\beta>0$ the process survives,
  then the process also survives for any $\beta'>\beta$.
\item If $\alpha$ and $\beta$ are fixed, and for some $\delta>0$ the process survives,
  then the process also survives for any $\delta'>\delta$.
\end{enumerate}
These properties follow easily from standard coupling arguments. We will denote by
$\beta_c=\beta_c(\alpha,\delta)\in[0,\infty]$ the parameter value such that, fixing these
$\alpha$ and $\delta$, $\eta_t$ survives for $\beta>\beta_c$ and dies out for
$\beta<\beta_c$. We define $\delta_c=\delta_c(\alpha,\beta)$ analogously.

Our first result provides some bounds on the critical parameters for survival. Let
$\beta^\text{\text{cp}}_c$ be the critical value 
of the contact process in $\zz^d$ (here we are taking the birth rate $\beta$ to be the
total birth rate from each site, so each site sends births to each given neighbor at rate
$\beta/(2d)$).

\begin{teo}\label{thm:parameters}\mbox{}
  \begin{enumerate}
  \item If $\beta\leq(\alpha+1)\beta^\text{\emph{cp}}_c$ then the process dies out.
  \item There exists a $\delta_p>0$ such that for any $\delta<\delta_p$ the process dies
    out (regardless of $\alpha$ and $\beta$).
  \item If $\beta>(\alpha+1)\beta^\text{\emph{cp}}_c$, then the process survives for large enough $\delta$.
  \end{enumerate}
\end{teo}

\begin{rem}
The published version of this article contains a different, more explicit version of (c), but it was pointed out to us by Dong Yao that its proof is flawed (the graphical representation employed in the construction of a contact process bounding our process from below does not enjoy all the necessary independence), and in fact the stated result appears to be wrong.
In this version we have thus replaced the statement with a weaker version of it, which however still exhibits an non-trivial phase transition between extinction and survival for our process.
The newly stated result can be proved analogously to Theorem 2.3 of \cite{LR}.
\end{rem}

Part (a) of the theorem is trivial, because the 1's die at rate $\alpha+1$. For part (b),
observe that if the complement of the set of sites at state $-1$ does not space-time
percolate, then each 1 in the process is doomed to live in a finite space-time region, and
then the process cannot have 1's at all times when started with finitely many occupied
sites. We will show by adapting arguments in \cite{meesterRoy} that, with probability 1,
no such space-time percolation occurs if $\delta$ is small enough.

A significant difficulty in giving a more complete picture of the phase diagram of
$\eta_t$ is that we lack a result about monotonicity with respect to $\alpha$ analogous to
the properties (i) and (ii) (monotonicity with respect to $\beta$ and $\delta$) mentioned
above. Observe that the equilibrium density of non-blocked sites is independent of
$\alpha$, but the environment changes more quickly as $\alpha$ increases. Simulations
suggest that if $\beta$ and $\delta$ are given and the process dies out at some parameter
value $\alpha$, then it also dies out for any parameter value $\alpha'>\alpha$ (note that
part (a) of Theorem \ref{thm:parameters} says that the process dies out at least for all
$\alpha$ large enough). But the usual simple arguments based on coupling do not work in
this case, since increasing $\alpha$ increases both the rate at which sites are blocked,
which plays against survival, and the rate at which sites are unblocked, which plays in
favor of survival, and we have not been able to find an alternative proof.

The second part of our study of $\eta_t$ investigates the convergence of the process and
the structure of its limit distributions. For $\eta\in\cx$ we will write $\eta=(A,B)$,
where
\[A=\{x\in\zz^d\!:\,\eta(x)=1\},\quad\text{and}\quad B=\{x\in\zz^d\!:\,\eta(x)=-1\}.\]
$\eta^\mu_t=(A^\mu_t,B^\mu_t)$ will denote the process with initial distribution $\mu$,
and we will refer to $B^\mu_t$ (or $B_t$ if no initial distribution is prescribed) as the
\emph{environment process}. Observe that the dynamics of the environment process are
independent of the 1's in $\eta_t$.

\begin{teo}\label{thm:compConv}
  Denote by $\tau=\inf\{t\geq0\!:\,A_t=\emptyset\}$ the extinction time of the
  process. Then for every initial distribution $\mu$,
  \[\eta^{\mu}_t\Longrightarrow\pp^{\mu}\!\left(\tau<\infty\right)\unu+\pp^{\mu}\!\left(\tau=\infty\right)\onu,\]
  where the limit is in the topology of weak convergence of probability measures.
\end{teo}

This result, which is usually called a \emph{complete convergence theorem}, implies that
all limit distributions are convex combinations of $\unu$ and $\onu$. Thus, the only
interesting non-trivial stationary distribution is $\onu$.

The proof of Theorem \ref{thm:compConv} relies on extending for $\eta_t$ the classical
block construction for the contact process introduced in \citet{bezGrimm}, so that we are
able to use the proof of complete convergence for the contact process to prove the
corresponding convergence of the contact process part of $\eta_t$. As a consequence of
this construction we will obtain, just as for the contact process, the fact that the
process dies out at the critical parameters $\beta_c$ and $\delta_c$ (see Corollary
\ref{cor:critical}). The arguments involved in this part will depend heavily on a duality
relation which will be developed in Section \ref{subsec:duality}.

The rest of the paper is devoted to the proofs of the two theorems.  Section
\ref{sec:preliminaries} describes the construction of $\eta_t$ and presents some basic
preliminary results. Theorem \ref{thm:parameters} is proved in Section \ref{sec:survExt}.
In Section \ref{sec:block} we obtain the block conditions for the survival of the process.
Finally, in Section \ref{sec:compConv} we use duality and the conditions obtained in
Section \ref{sec:block} to prove Theorem \ref{thm:compConv}.

\section{Preliminaries}\label{sec:preliminaries}
\renewcommand{\theteo}{\arabic{teo}}

\subsection{Graphical representation and monotonicity}\label{subsec:constMon}

The graphical representation is one of the basic and most useful tools in the study of the
contact process and other interacting particle systems. It will allow us to construct our
process from a collection of independent Poisson processes and obtain a single probability
space in which copies of the process with arbitrary initial states can be coupled.  We
will give a rather informal description of this construction, which can be made precise by
adapting the arguments of \citet{harris2}. We refer the reader to Section III.6 of
\citet{ligg1} for more details on this construction in the case of an additive spin
system.

The construction is done by placing symbols in $\zz^d\times[0,\infty)$ to represent the
different events in the process. For each ordered pair $x,y\in\zz^d$ at distance 1 let
$N^{x,y}$ be a Poisson process with rate $\beta/(2d)$, and take the processes
assigned to different pairs to be independent. At each event time $t$ of $N^{x,y}$ draw an
arrow $\xrightarrow[\quad]{1}$ in $\zz^d\!\times\![0,\infty)$ from $(x,t)$ to $(y,t)$ to
indicate the birth of a 1 sent from $x$ to $y$ (which will only take place if $x$ is
occupied and $y$ is vacant at time $t$). Similarly, define a family of independent Poisson
processes $(U^{1,x})_{x\in\zz^d}$ with rate 1 and for each event time $t$ of $U^{1,x}$
place a symbol $\ast_1$ at $(x,t)$ to indicate that a 1 flips to 0 (i.e. that a particle
dies). To represent the environment, consider two families of independent Poisson
processes $(V^{x})_{x\in\zz^d}$ and $(U^{-1,x})_{x\in\zz^d}$ with rates $\alpha$ and
$\alpha\delta$ respectively. For each event time $t$ of $V^x$ place a symbol
$\bullet_{-1}$ at $(x,t)$ to indicate the birth of a $-1$ (i.e. the blocking of a site)
and for each event time $t$ of $U^{-1,x}$ place a symbol $\ast_{-1}$ to indicate that a
$-1$ flips to 0 (i.e. the unblocking of a site).

We construct $\eta_t$ from this percolation structure in the following way. Consider a
deterministic initial condition $\eta_0$ and define the environment process $B_t$ by
setting $\eta_t(x)=-1$ when $(x,t)$ lies between symbols $\bullet_{-1}$ and $\ast_{-1}$
(in that order) in the time line $\{x\}\times[0,\infty)$, and also if $\eta_0(x)=-1$ and
there is no symbol $\ast_{-1}$ in that time line before time $t$. Having defined $B_t$, we
say that there is an \emph{active path} between $(x,s)$ and $(y,t)$ if there is a
connected oriented path, moving along the time lines in the increasing direction of time
and passing along arrows $\xrightarrow[\quad]{1}$, which crosses neither symbols $\ast_1$
nor space-time points that were set to $-1$. The collection of active paths corresponds to
the possible space-time paths along which 1's can move, so we define $A_t$ by
\[A_t=\{y\in\zz^d\!:\,\exists x\in A_0\text{ with an active path from $(x,0)$ to
  $(y,t)$}\}.\] The arguments of \citet{harris2} imply that this construction gives a well
defined Markov process with the right transition rates. Moreover, the same realization of
this graphical representation can be used for different initial conditions, and this gives
the coupling mentioned above (see Section III.6 in \citet{ligg1} for more details on this
coupling in the case of a spin system). For the rest of the paper we will implicitly use
this ``canonical'' coupling every time we couple copies of $\eta_t$ with different initial
conditions. The attractiveness property mentioned in the \nameref{sec:intr} follows
directly from this construction, and the monotonicity properties with respect to $\beta$
and $\delta$ can be obtained by a simple modification of this coupling (analogous to what
is done for the contact process).

Recall the definition of the partial order on configurations given in \eqref{eq:order}.
Clearly,
\[\eta^1\leq\eta^2\quad\Leftrightarrow\quad A^1\subseteq
A^2\text{ and }B^1\supseteq B^2.\] For probability measures on $\cx$, which we endow with
the product topology, we consider the usual ordering: $\mu_1\leq\mu_2$ if and only if
$\int\!f\,d\mu_1\leq\int\!f\,d\mu_2$ for every continuous increasing
$f\!:\cx\longrightarrow\rr$. We recall that the property $\mu_1\leq\mu_2$ is equivalent to
the existence of a probability space in which a pair of random variables $X_1$ and $X_2$
with distributions $\mu_1$ and $\mu_2$ can be coupled in such a way that $X_1\leq X_2$
almost surely (see Theorem II.2.4 in \citet{ligg1}). We will use this fact repeatedly, and
for simplicity we will say that $X_2$ \emph{dominates} $X_1$ when this condition holds. We
will also use this term to compare two processes, so saying that $\eta^2_t$ dominates
$\eta^1_t$ will mean that the two processes can be constructed in a single probability
space in such a way that $\eta^1_t\leq\eta^2_t$ for all $t\geq0$.

The attractiveness property allows us to obtain the lower and upper invariant measure of
the process.

\begin{prop}\label{prop:mon}
  Let $\chi_{\zz^d}$ be the probability distribution on $\cx$ assigning mass $1$ to the
  all $1$'s configuration, and let $S(t)$ be the semigroup associated to the process.
  Define
  \[\onu=\lim_{t\rightarrow\infty}\chi_{\zz^d}S(t),\]
  where the limit is in the topology of weak convergence of probability measures. Then
  $\onu$ is the upper invariant measure of the process, that is, $\onu$ is invariant and
  every other invariant measure is stochastically smaller than $\onu$. Moreover,
  \[\onu=\lim_{t\rightarrow\infty}\nu_{\zz^d}S(t).\]

  Analogously,
  \[\unu=\nu_\emptyset\]
  is the lower invariant measure of the process.
\end{prop}

\begin{proof}
  Since $\mr$ is invariant for the environment and the empty state is a trap for the 1's,
  $\unu$ is invariant. It is the lower invariant measure because every invariant measure
  has $\mr$ as its projection onto the environment, and $\nu_\emptyset$ is the smallest
  probability measure on $\cx$ having $\mr$ as its marginal on the $-1$'s.

  For the other part, standard arguments imply that the limit defining $\onu$ exists and
  is invariant (see, for instance, Sections I.1 and III.2 in \citet{ligg1}). Since
  $\chi_{\zz^d}$ is larger than any other measure on $\cx$, it follows by attractiveness
  that $\onu$ is the largest invariant measure.

  Now let $\nu^*=\lim_{t\rightarrow\infty}\nu_{\zz^d}S(t)$. As above, $\nu^*$ is
  well-defined and invariant, so to prove that $\nu^*=\onu$ it is enough to prove that
  $\nu^*$ is larger than any other invariant measure. If $\nu$ is any invariant measure,
  its projection onto the $-1$'s must be $\mr$, so for any continuous increasing $f$,
  \begin{align*}
    \int\!f\,d\nu=\ee^\nu\!\left(f(\eta_0)\right)&=\ee^\nu\!\left(f(\eta_t)\right)\\
    &\leq\ee^{\nu_{\zz^d}}\!\left(f(\eta_t)\right)=\int\!f\,d[\nu_{\zz^d}S(t)]
    \xrightarrow[t\rightarrow\infty]{}\int\!f\,d\nu^*.\qedhere
  \end{align*}
\end{proof}

\subsection{Duality}\label{subsec:duality}

The dual process $(\dual{\eta}{t}_s)_{0\leq s\leq t}=(\dual{A}{t}_s,\dual{B}{t}_s)_{0\leq
  s\leq t}$ is constructed using the same graphical representation we used for
constructing $\eta_t$. Our duality relation will require that the process be started with
the environment at equilibrium. The dual processes will also be started with measures of
the form $\nC$, for $C\subseteq\zz^d$, and the dual process started with this distribution
will be denoted by $(\dual{\eta}{\nC,t}_s)_{0\leq s\leq t}$.

Fix $t>0$, and start by choosing $B_0$ according to $\mr$. Then run the environment
process forward in time until $t$, using the graphical representation. This defines
$(B_s)_{0\leq s\leq t}$. The dual environment is given by $\dual{B}{t}_s=B_{t-s}$. Now
place a 1 at time $t$ at every site $x\in C\setminus \dual{B}{t}_0$, that is, every site
in $C$ which is not blocked by the environment at time $t$. This defines
$\dual{A}{\nC,t}_0$, and by the stationarity of the environment process we get an initial
condition $(\dual{A}{\nC,t}_0,\dual{B}{t}_0)$ for the dual chosen according to $\nC$.
Having defined $\dual{A}{\nC,t}_0$ and $(\dual{B}{t}_s)_{0\leq s\leq t}$, we define the
1-dual by
\[\dual{A}{\nC,t}_s=\{y\in\zz^d\!:\,\exists x\in\dual{A}{\nC,t}_0\text{ with an active
  path from $(y,t-s)$ to $(x,t)$}\}.\] That is, the 1-dual is defined by running the
contact process for the 1's backwards in time and with the direction of the arrows
reversed. An active path in $\eta_t$ from $(y,t-s)$ to $(x,t)$ will be called a \emph{dual
  active path} from $(x,t)$ to $(y,t-s)$ in the dual process.

We could have defined the dual by simply choosing a random configuration at time $t$
according to $\nC$ and then running the whole process backwards. The idea of the preceding
construction is to allow coupling the process and its dual in the same graphical
representation in such a way that the initial state of the environment for $\eta_s$ is the
same as the final state of the environment for $\dual{\eta}{t}_s$ (that is,
$B_0=\dual{B}{t}_t$). This allows us to obtain the following duality result:

\begin{prop}\label{prop:dual}
  For any $A,C,D\subseteq\zz^d$,
  \begin{equation}\label{eq:dual1}
    \pp^\nA\!\left(\vphantom{\dual{A}{t}}A_t\cap C\notempty,B_t\cap
      D\notempty\right)=\pp^\nC\!\left(\dual{A}{t}_t\cap
      A\notempty,\dual{B}{t}_0\cap D\notempty\right).
  \end{equation}
  Moreover, $\eta_t$ satisfies the following self-duality relation: if $A$ or $C$ is
  finite, then
  \begin{equation}\label{eq:dual2}
    \pp^\nA\!\left(A_t\cap C\notempty,B_t\cap
      D\notempty\right)=\pp^\nC\!\left(A_t\cap A\notempty,B_0\cap
      D\notempty\right).
  \end{equation}
\end{prop}

\begin{proof}
  The first equality follows directly from coupling the process and its dual using the
  same realization of the graphical representation.  Indeed, if we use this coupling then, by
  definition,
  \[\pp\!\left(\dual{B}{\nC,t}_s=B^\nA_{t-s}\text{ for every $0\leq s\leq
      t$}\right)=1.\] Calling $\cal E$ the $\sigma$-algebra generated by the environment
  process, observe that our construction implies that
  \[\pp^\nA(A_t\cap C\notempty|{\cal E})=\pp^\nC(\dual{A}{t}_t\cap A\notempty|{\cal E}).\]
  Therefore,
  \begin{align*}
    \pp^\nA\!\left(\vphantom{\dual{A}{t}}A_t\cap C\notempty,B_t\cap
      D\notempty\right)&=\ee^\nA\!\left(\vphantom{\dual{A}{t}}\pp(A_t\cap C\notempty|{\cal
        E}),B_t\cap
      D\notempty\right)\\
    &=\ee^\nC\!\left(\pp(\dual{A}{t}_t\cap A\notempty|{\cal E}),\dual{B}{t}_0\cap D\notempty\right)\\
    &=\pp^\nC\!\left(\dual{A}{t}_t\cap A\notempty,\dual{B}{t}_0\cap D\notempty\right).
  \end{align*}

  \eqref{eq:dual2} is obtained from \eqref{eq:dual1}, the self-duality of the contact
  process, and the reversibility of the environment.
\end{proof}

Taking $A$ finite and $C=D=\zz^d$ in \eqref{eq:dual2} and using the monotonicity of the
event $\{A_t\notempty\}$ in $t$ we obtain the following:
\begin{equation*}
  \pp^\nA\!\left(A_t\notempty\,\,\forall
    t\geq0\right)=\onu\left(\{(E,F)\!:\,E\cap A\notempty\}\right).
\end{equation*}
Since $\onu$ is translation invariant, the right side of this equality is positive if and
only if $A\notempty$ and $\eta_t$ survives, that is, $\onu\neq\unu$. As a consequence we
deduce that the following condition is equivalent to the survival of the process:
\begin{equation}\tag{S1}\label{S1}
  \parbox{4.5in}{For any (or, equivalently, some) finite $A\subseteq\zz^d$ with
    $A\notempty$, the process started at $\nA$ contains 1's for every $t\geq0$ with positive probability.}
\end{equation}

\subsection{Positive correlations}

A second property that is central to the study of the contact process is \emph{positive
  correlations}. Recall that a probability measure $\mu$ has positive correlations if for
every $f,g$ increasing,
\begin{equation}\label{eq:posCorr}
  \int\!fg\,d\mu\geq\int\!f\,d\mu\int\!g\,d\mu.
\end{equation}

In the following lemma we prove a version of positive correlations for $\eta^\nA_t$ with
respect to cylinder functions.

\begin{lem}\label{lem:posCorr}
  Let $f,g$ be increasing real-valued functions on $\cx$ depending on finitely many
  coordinates. Then if $\mu_t$ denotes the distribution of $\eta_t^\nA$,
  \eqref{eq:posCorr} holds with $\mu=\mu_t$, that is,
  \begin{equation}\label{eq:finPosCorr}
    \ee^\nA\!\left(f(\eta_t)g(\eta_t)\right)\geq\ee^\nA\!\left(f(\eta_t)\right)\ee^\nA\!\left(g(\eta_t)\right).
  \end{equation}
  The same inequality holds if $\nA$ is replaced by any deterministic initial condition.
\end{lem}

\begin{proof}
  Since $f$ and $g$ depend on finitely many coordinates and every jump in our process is
  between states which are comparable in the partial order \eqref{eq:order}, a result of
  Harris (see Theorem II.2.14 in \citet{ligg1}) and attractiveness imply that it is enough
  to show that the initial distribution of the process has positive correlations in the
  sense of the lemma. The result with $\nA$ replaced by a deterministic initial condition
  readily follows.

  To show that $\nA$ is positively correlated, consider the process $\varsigma_t$ defined
  in $\cx$ by $\varsigma_0\equiv1$ and independent transitions at each site given by
  \begin{align*}
    \left.\begin{aligned}
        &0\longrightarrow & -1  &\quad\text{at rate $\rho$}\\
        -&1\longrightarrow & 0 &\quad\text{at rate $1-\rho$}
      \end{aligned}\,\right\}\,&\text{for $x\notin A$,}\\
    \left.\begin{aligned}
        &1\longrightarrow & -1  &\quad\text{at rate $\rho$}\\
        -&1\longrightarrow & 1 &\quad\text{at rate $1-\rho$}
      \end{aligned}\,\right\}\,&\text{for $x\in A$.}
  \end{align*}
  It is clear that $\varsigma_t$ converges weakly to the measure $\nA$. Since the initial
  distribution of $\varsigma_t$ has positive correlations (because it is deterministic),
  \eqref{eq:posCorr} holds for its limit $\nA$, using again Harris' result.
\end{proof}

\section{Survival and extinction}\label{sec:survExt}

In this section we prove Theorem \ref{thm:parameters}. Throughout the proof we will
implicitly use \eqref{S1} to characterize survival. We start with the easy part.

\begin{proof}[Proof of Theorem \ref{thm:parameters}, part (a)]
  Consider the process $\widetilde{\eta}_t$ defined by the following transition rates:
  \begin{center}
    \begin{tabular}{ccclc}
      $0,-1$ & $\longrightarrow$ & $1$\hspace{.3in} & at rate & $\beta f_1$ \\
      $1$ & $\longrightarrow$ & $0$\hspace{.3in} & at rate & 1 \\
      $0,1$ & $\longrightarrow$ & $-1$\hspace{.3in} & at rate & $\alpha$ \\
      $-1$ & $\longrightarrow$ & $0$\hspace{.3in} & at rate & $\alpha\delta$ \\
    \end{tabular}
  \end{center}
  This process corresponds to modifying $\eta_t$ by ignoring the effect of blocked sites
  on births. It is easy to couple $\widetilde{\eta}_t$ and $\eta_t$ using the graphical
  representation in such a way that if the initial states are the same,
  $\eta_t\leq\widetilde{\eta}_t$ for all $t\geq0$. Therefore, it is enough to show that
  $\widetilde{\eta}_t$ dies out, and this follows directly from the hypothesis because the
  1's in $\widetilde{\eta}_t$ behave just like a contact process with birth rate $\beta$
  and death rate $\alpha+1$.
\end{proof}

The proof of part (b) is more involved, and it is based on adapting the techniques of
Boolean models in continuum percolation (see \citet{meesterRoy}).

\begin{proof}[Proof of Theorem \ref{thm:parameters}, part (b)]
  The idea is to show that when $\delta$ is small, the set of unblocked sites in the
  environment process $B_t$ does not ``space-time percolate'' with probability 1. By this
  we mean that there is no infinite path in $\zz^d\times[0,\infty)$ moving between
  nearest-neighbor sites in $\zz^d$ and along time lines in the increasing direction of
  time that uses only non-blocked sites. The conclusion follows directly from this fact,
  since in that case every 1 will live in a finite space-time box, so it will not be able
  to contribute to the survival of the process.

  By a simple time change, we can consider the environment process as having transitions
  given by
  \begin{center}
    \begin{tabular}{ccclc}
      $-1$ & $\longrightarrow$ & $0$\hspace{.3in} & at rate & $q$ \\
      $0$ & $\longrightarrow$ & $-1$\hspace{.3in} & at rate & $1-q$,
    \end{tabular}
  \end{center}
  where $q=\delta/(1+\delta)\longrightarrow0$ as $\delta\longrightarrow0$.  We still
  consider this process as defined by the graphical representation, though now the symbols
  $\bullet_{-1}$ and $\ast_{-1}$ appear at rate $1-q$ and $q$ respectively.

  Take the percolation structure given by the graphical representation and draw for every
  symbol $\ast_{-1}$ at a space-time point $(x,t)$ a box of base $x+[-2/3,2/3]^d$ spanning
  the interval in the time coordinate from $t$ until the time corresponding to the next
  symbol $\bullet_{-1}$ (i.e., these boxes span intervals where the sites are not
  blocked). Then, since the environment process is translation invariant, the 0's will
  almost surely not space-time percolate if and only if
  \begin{equation}\label{eq:noperc}
    \pp(|\mathcal{W}|=\infty)=0,
  \end{equation}
  where $\mathcal{W}$ denotes the connected component of the union of the boxes that
  contains the origin at time 0, and $|\mathcal{W}|$ denotes the number of boxes that form
  this cluster.

  To prove \eqref{eq:noperc} we compare this continuum percolation structure with a
  multitype branching process $X=(X_{n,i})_{n,i\in\nn}$. The first step in the comparison
  is to stretch all the boxes so that their heights are all integer-valued. It is enough
  to show that \eqref{eq:noperc} holds after this modification, since increasing the
  heights of the boxes increases the probability of space-time percolation of the
  unblocked sites. Assume that the origin is not blocked at time 0, and call $i_0\in\nn$ the
  (random) height of its associated box. For simplicity, assume further that all the
  neighbors of the origin are blocked at time 0, the extension to the general case being
  straightforward. We start defining $X$ by saying that the 0-th generation has only one
  member, and it is of type $i_0$ (that is, $X_{0,j}=\uno{\{j=i_0\}}$). The box containing
  the origin at time 0 is possibly intersected by boxes placed at the $2d$ neighbors of
  the origin, and these boxes will constitute the children of the initial member: we let
  $X_{1,j}$ be the number of boxes of height $j$ that intersect the original box.  We
  define the subsequent generations of $X$ inductively: $X_{n+1,j}$ is the number of boxes
  of height $j$ that intersect boxes of the $n$-th generation and which have not been
  counted up to generation $n-1$.  Now let
  \[X^\infty=\sum_{n=0}^\infty\sum_{i=1}^\infty X_{n,i},\] and observe that every box in
  $\mathcal{W}$ is counted in $X^\infty$, so
  \begin{equation}\label{eq:bdCardW}
    |\mathcal{W}|\leq X^\infty
  \end{equation}
  (recall that $X$ is constructed from the stretched boxes).
  
  Our goal is to show that $\ee(X^\infty)<\infty$. To achieve this we will couple $X$ with
  another multitype branching process $Y=(Y_{n,i})_{n,i\in\nn}$, which we define
  below. The details of this part can be adapted easily from the proof of Theorem 3.2 in
  \citet{meesterRoy}, so we will only sketch the main ideas. Consider a box of height $i$
  based at $[x-2/3,x+2/3]^d\times\{t\}$, which we will denote by $B(x,t,i)$. The boxes of
  height $j$ that intersect this box must all have bases of the form
  $[y-2/3,y+2/3]\times\{s\}$ for some $y$ at distance 1 of $x$ and some
  $s\in(0\vee(t-j),t+i]$. The number of symbols $\ast_{-1}$ appearing in the piece
  $\{y\}\times(0\vee(t-j),t+i]$ of the graphical representation above a given neighbor $y$
  of $x$ is a Poisson random variable with mean $q[t+i-0\vee(t-j)]\leq q[i+j]$, and
  each of these symbols corresponds to a box that intersects $B(x,t,i)$. Since the
  probability that any one of these (stretched) boxes is of height $j$ is
  $p_j=\pp\left(Z\in(j-1,j]\right)$, where $Z$ is an exponential random variable with rate
  $1-q$, we deduce that the number of children of $B(x,t,i)$ of height $j$ is a Poisson
  random variable with mean bounded by
  \begin{equation}\label{eq:bdRateBox}
    2dp_jq[i+j]\leq4dqijp_j,
  \end{equation}
  where we used the fact that $i+j\leq2ij$ for positive integers $i$ and $j$. Now let $Y$
  be a multitype branching process where the number of children of type $j$ of each
  individual of type $i$ is a Poisson random variable with mean $4dqijp_j$ ($Y_{n,i}$ is
  the number of individuals of type $i$ in generation $n$). Then a coupling argument and
  \eqref{eq:bdRateBox} imply that if $X_{0,i}=Y_{0,i}$ for all $i\geq1$ then $X_{n,i}$ is
  dominated by $Y_{n,i}$ for each $n\geq0$ and $i\geq1$, and thus
  \begin{equation}\label{eq:bdXInfty}
    \ee\Big(X^\infty\big|X_{0,k}=\uno{\{k=i_0\}}\Big)
    \leq\ee\bigg(\sum_{n=0}^\infty\sum_{j=1}^\infty
    Y_{n,j}\Big|Y_{0,k}=\uno{\{k=i_0\}}\bigg).
  \end{equation}
  To bound this last sum we recall a standard result in branching processes theory (see,
  for example, Chapter V in \citet{athreyaNey}): the expected number of individuals of
  type $j$ in the $n$-th generation of $Y$ when starting with one individual of type $i_0$
  is given by
  \begin{equation}\label{eq:YnMn}
    \ee\Big(Y_{n,j}\big|Y_{0,k}=\uno{\{k=i_0\}}\Big)=(M^n)_{i_0,j},
  \end{equation}
  where $M$ is the infinite matrix indexed by $\nn$ with $M_{i,j}$ being
  the expected number of children of type $j$ of an individual of type $i$. By definition
  of $Y$, $M_{i,j}=4dqijp_j$, and from this we get inductively a bound for
  $(M^n)_{i_0,j}$:
  \[(M^n)_{i_0,j}\leq(4dq)^ni_0\ee\big(H^2\big)^{n-1}\pp(H=j)j\] for all $n\geq1$, where
  $H$ is a random variable with positive integer values and distribution given by
  $\pp(H=j)=p_j$. Using this together with \eqref{eq:bdXInfty} and \eqref{eq:YnMn} gives
  \begin{equation}\label{eq:boundX}
    \begin{aligned}
      \ee\left(\vphantom{2^2}X^\infty\right|
      \left.\vphantom{2^2}X_{0,k}=\uno{k=i_0}\right)&\leq
      1+i_0\sum_{n=1}^\infty\!\left((4dq)^n\ee\big(H^2\big)^{n-1}\sum_{j=1}^\infty
        p_jj\right)\\
      &=1+4dqi_0\ee(H)\sum_{n=0}^\infty(4dq\ee(H^2))^n.
    \end{aligned}
  \end{equation}
  Observe that $H$ is dominated by $Z+1$, so
  $\ee\big(H^2\big)\leq\frac{2(2-q)}{(1-q)^2}+1$.  Hence,
  \begin{equation}\label{eq:boundDp}
    4dq\ee(H^2)\leq4d\left(\frac{2q(2-q)}{(1-q)^2}+q\right)<1
  \end{equation}
  for sufficiently small $q$, and then the last sum in \eqref{eq:boundX} converges for such
  $q$. This implies by \eqref{eq:bdCardW} that $\ee(|\mathcal{W}|)<\infty$, so
  $\pp(|\mathcal{W}|=\infty)=0$.
\end{proof}

Using \eqref{eq:boundDp} we can get explicit lower bounds for $\delta_p$, but these turn
out to be rather small (around $0.02$ for $d=2$ and $0.01$ for $d=3$).

\begin{proof}[Proof of Theorem \ref{thm:parameters}, part (c)]
  The effect on our process of letting $\delta\to\infty$ is that blocked sites become immediately unblocked, and thus blocking events have just the same effect as death of the particle possibly occupying the site.
  This heuristic suggests that as $\delta\to\infty$ our process behaves simply as the usual contact process with birth rate $\beta$ and death rate $\alpha+1$, which implies the claimed result.
  This idea can be turned into a proof roughly as follows. If we look at a fixed space time box, then with high probability (as $\delta\to\infty$ with all other parameters being fixed) the following holds: every site that gets blocked inside the block has to become unblocked before any infection event enters of exits the site.
  Under this condition the process never really sees any blocked sites, and thus it behaves as the desired usual contact process.
  One can then use a block construction as usual to show survival, asking that the ``good blocks'' in this argument satisfy the specified condition.
  Implementing this idea rigorously takes a bit of additional work, but it can be done very similarly to the proof of Theorem 2.3 in \cite{LR}, so we omit the details.
\end{proof}

\section{Block construction}\label{sec:block}

The aim of this section is to establish ``block conditions'' concerning the process in a
finite space-time box that guarantee survival. This was first done in \citet{bezGrimm}.
Here we will follow closely Section I.2 of \citet{ligg2}, together with the corrections to
the book that can be found in the author's website.

Before getting started with the block construction we need to obtain the equivalent
condition for survival mentioned in the \nameref{sec:intr}, which says that $\eta_t$
survives if and only if the following condition holds:
\begin{equation}\tag{S2}\label{S2}
  \parbox{4.5in}{The process started with a single 1 at the origin and everything else at
    $-1$ contains 1's at all times with positive probability.}
\end{equation}
The sufficiency of this condition is a consequence of \eqref{S1} and attractiveness. The
necessity will be a consequence of the following stronger result, which is precisely what
we will need in the proof of Lemma \ref{lem:techBlock} below. Let $\chi_A$ denote the
probability measure on $\cx$ that assigns mass 1 to the configuration $\eta$ with
$\eta|_A\equiv1$, $\eta|_{A^c}\equiv-1$.

\begin{lem}\label{lem:bigN}
  Suppose that the process survives. Then for any $\sigma>0$ there is a positive integer
  $n$ such that
  \begin{equation*}
    \pp^{\chi_{[-n,n]^d}}\!\left(A_t\notempty\,\,\forall
      t\geq0\right)>1-\sigma^2.
  \end{equation*}
\end{lem}

To obtain \eqref{S2} from this result observe that the process started with a single 1 at
the origin has $[-n,n]^d$ fully occupied by time 1 with some positive probability, so we
can use the strong Markov property and attractiveness to restart the process at time 1
starting from $\chi_{[-n,n]^d}$ and obtain
$\pp^{\chi_{\{0\}}}\!\left(A_t\notempty\,\,\forall t\geq0\right)>0$.  Observe that the
lemma is a simple consequence of duality when the initial condition for $\eta_t$ is
$\nu_{[-n,n]^d}$ instead of $\chi_{[-n,n]^d}$.  Indeed, using \eqref{eq:dual2} with
$D=\zz^d$ gives
\begin{align*}
  \hspace{0.4in}&\hspace{-0.4in}\lim_{n\rightarrow\infty}\pp^{\nu_{[-n,n]^d}}\!\left(A_t\notempty\,\,\forall
    t\geq0\right)\\
  &=\lim_{n\rightarrow\infty}\lim_{t\rightarrow\infty}\pp^{\nu_{[-n,n]^d}}\!\left(A_t\notempty\right)
  =\lim_{n\rightarrow\infty}\lim_{t\rightarrow\infty}\pp^{\nu_{\zz^d}}\!\left(A_t\cap[-n,n]^d\notempty\right)\\
  &=\lim_{n\rightarrow\infty}\onu\left(\vphantom{2^2}\{(E,F)\!:\,E\cap[-n,n]^d\notempty\}\right)
  =\onu\left(\vphantom{2^2}\{(E,F)\!:\,E\notempty\}\right).
\end{align*}
This last probability is 1 when $\eta_t$ survives, so in this case given any
$\varepsilon>0$ we can choose a positive integer $m$ such that
\begin{equation}\label{eq:bigM}
  \pp^{\nu_{[-m,m]^d}}\!\left(A_t\notempty\,\,\forall
    t\geq0\right)>1-\varepsilon.
\end{equation}

Recall that in Proposition \ref{prop:mon} we showed that the limit distributions of the
processes started at $\chi_{\zz^d}$ and at $\nu_{\zz^d}$ are the same. It is then
reasonable to expect that the asymptotic behavior as $t\rightarrow\infty$ of the process
started at $\chi_{[-n,n]^d}$ is similar to that of the process started at
$\nu_{[-n,n]^d}$, at least for large enough $n$. This idea will allow us to derive the
lemma from \eqref{eq:bigM}.

\begin{proof}[Proof of Lemma \ref{lem:bigN}]
  Let $\varepsilon>0$ and choose $m$ to be the positive integer obtained in
  \eqref{eq:bigM}. To extend this inequality to the process started at $\chi_{[-n,n]^d}$
  we will consider two copies of the process $\eta^1_t$ and $\eta^2_t$ coupled using the
  graphical representation, with $\eta^1_t$ started at $\nu_{[-m,m]^d}$ and $\eta^2_t$ at
  $\chi_{[-n,n]^d}$ for some large $n>m$. For simplicity we will write $Q(k)=[-k,k]^d$.

  We want to obtain a space-time cone growing linearly in time such that
  $\cup_{t\geq0}\{t\}\!\times\!A^{\nu_{Q(m)}}_t$ is contained in that cone with high probability. To
  achieve this we compare $A^{\nu_{Q(m)}}_t$ with a branching random walk $Z_t$ with
  branching rate $\beta/(2d)$ and no deaths (that is, each particle in $Z_t$ gives birth
  to a new particle at each neighbor at rate $\beta/(2d)$, and multiple particles per site
  are allowed). Let $\{p_t(x,y)\}_{x,y\in\zz^d}$ be the transition probabilities of a
  simple random walk in $\zz^d$ that moves to each neighbor at rate $\beta/(2d)$ and let
  $C_t$ be the set-valued process given by
  \[C_t=\big\{x\in\zz^d\!:\,Z_t(x)>0\big\}.\] For $D\subseteq\zz^d$, $Z^D_t$ and $C^D_t$
  will denote the processes started with all sites in $D$ occupied by one particle and no
  particles outside $D$. It is not hard to see that for any $t>0$ and any $x\in\zz^d$,
  \[\ee\!\left(Z^{\{0\}}_t(x)\right)=e^{\beta t}p_t(0,x)\]
  (see, for instance, the proof of Proposition I.1.21 in \citet{ligg2}). Therefore, for
  any $D\subseteq\zz^d$,
  \[\ee\!\left(\left|C^{\{0\}}_t\cap D^\text{c}\right|\right)
  \leq \sum_{x\notin D}\ee\!\left(Z^{\{0\}}_t(x)\right) =e^{\beta t}\sum_{x\notin
    D}p_t(0,x).\] From this we get that if $k>m$ and $c>0$ then
  \begin{equation}\label{eq:exCube}
    \ee\!\left(\left|C^{Q(m)}_t\cap Q(k+ct)^\text{c}\right|\right)
    \leq(2m+1)^de^{\beta t}\!\sum_{\|x\|_\infty>k-m+ct}p_t(0,x).
  \end{equation}
  Now if $X_t$ is the one dimensional random walk starting at $0$ and moving to each
  neighbor at rate $\beta/(2d)$, Chebyshev's inequality gives
  \begin{align*}
    \pp(|X_t|>k-m+ct)&=2\pp(X_t-k+m-ct>0)\\
    &\leq2\ee\!\left(e^{X_t-k+m-ct}\right)=2e^{-(k-m)}e^{-ct+\frac{\beta}{2d}(e+e^{-1}-2)t}.
  \end{align*}
  The last equality can be obtained by seeing $X_t$ as the difference between two
  independent Poisson random variables, each with mean $(\beta t)/(2d)$, and using the
  fact that the moment generating function of a Poisson random variable $Y$ with mean
  $\lambda$ is $\ee\!\left(e^{sY}\right)=e^{\lambda(e^{s}-1)}$.  Applying this bound to
  each coordinate of the $d$-dimensional walk we get that
  \[\sum_{\|x\|_\infty>k-m+ct}p_t(0,x)\leq d\,\pp(|X_t|>k-m+ct)\leq
  2de^{-(k-m)}e^{-ct+\frac{\beta}{2d}(e+e^{-1}-2)t},\] and
  then using \eqref{eq:exCube} we deduce that $c$ can be taken large enough so that
  \begin{equation*}
    \ee\!\left(\left|C^{Q(m)}_t\cap Q(k+ct)^\text{c}\right|\right)
    \leq 2d(2m+1)^de^{-(k-m)}e^{-t}.
  \end{equation*}
  Observe that, by the definition of $Z_t$, the process $A^{\nu_{Q(m)}}_t$ is dominated by
  $C^{Q(m)}_t$, so the last bound implies that
  \begin{multline}\label{eq:intExpec}
    \ee\Bigg(\int_0^\infty\!\big|A_t^{\nu_{Q(m)}}\cap Q(k+ct)^\text{c}\big|\,dt\Bigg)\\
    \leq\int_0^\infty\!\ee\!\left(\big|C^{Q(m)}_t\cap Q(k+ct)^\text{c}\big|\right)dt
    \leq2d(2m+1)^de^{-(k-m)}.
  \end{multline}
  We can use this inequality to estimate the probability that $A_t\subseteq Q(k+1+ct)$ for
  all $t\geq 0$. Observe that if $x\in A_t\cap Q(k+1+ct)^\text{c}$, the particle at $x$
  survives at least until time $t+2/c$ with probability $e^{-2\alpha(1+\delta)/c}$, and
  thus $x\in A_s\cap Q(k+cs)^\text{c}$ for all $s\in[t+1/c,t+2/c]$ with at least that
  probability. We deduce that
  \begin{multline*}
    \ee^{\nu_{Q(m)}}\!\left(\int_0^\infty\!\left|A_t\cap Q(k+ct))^\text{c}\right|dt\right)\\
    \geq\pp^{\nu_{Q(m)}}\!\left(A_t\cap Q(k+1+ct)^\text{c}\notempty\,\,\text{for some
      }t\geq0\right)e^{-2\alpha(1+\delta)/c}\frac{1}{c}.
  \end{multline*}
  Therefore, if we let
  \[\mathrm{G}_1=\left\{A^1_t\subseteq Q(k+1+ct)\,\,\forall t\geq0\right\},\]
  (where $A^1_t$ denotes the set of 1's in the process $\eta^1_t$ started at
  $\nu_{Q(m)}$), we can use this bound together with \eqref{eq:intExpec} to get
  \[\pp(\mathrm{G}_1^\text{c})\leq 2cd(2m+1)^de^{2\alpha(1+\delta)/c}e^{-(k-m)}.\]
  Choosing now $k$ large enough yields
  \begin{equation}\label{eq:cone}
    \pp\!\left(\mathrm{G}_1\right)>1-\varepsilon.
  \end{equation}

  Now take $n>k$, $T>0$, let $(t-T)^+=(t-T)\vee0$, and call $\mathrm{G}_2$ the event that
  on the space-time region $\cup_{t\geq0}\{t\}\!\times\!Q(n+c(t-T)^+)$ the environment for $\eta^2_t$
  dominates the environment for $\eta^1_t$ (with respect to the order \eqref{eq:order}):
  \[\mathrm{G}_2=\left\{B^2_t\subseteq B^1_t\,\,\text{on
    }Q(n+c(t-T)^+)\,\,\forall t\geq0\right\}.\] We want this space-time region to contain
  the region defining $\mathrm{G}_1$, so we let $T=(n-k-1)/c$.

  Observe that, since we are coupling the processes using the canonical coupling given by
  the graphical representation, once the environment is equal for both process at a given
  site, it stays equal at that site from that time on. In particular, $B^2_t$ dominates
  $B^1_t$ on $Q(n)$ for all $t\geq0$.  For any other site, any symbol $\bullet_{-1}$ or
  $\ast_{-1}$ leaves the environment equal for both process. Therefore,
  \begin{align*}
    \pp(\mathrm{G}_2^\text{c})&\leq\sum_{x\notin Q(n)}\pp\bigg(\text{no $\bullet_{-1}$ or
      $\ast_{-1}$
      at $x$ by time }T+(\|x\|_\infty-n)/c\bigg)\\
    &=\sum_{j>n}|Q(j)\setminus
    Q(j-1)|e^{-\alpha(1+\delta)(T+(j-n)/c)}\\
    &\leq e^{\alpha(1+\delta)(k+1)/c}\sum_{j>n}(2j+1)^d e^{-\alpha(1+\delta)j/c}.
  \end{align*}
  By taking $n$ large enough we obtain
  \begin{equation}\label{eq:cone2}
    \pp(\mathrm{G}_2)>1-\varepsilon.
  \end{equation}

  Finally, let
  \[\mathrm{G}_3=\left\{A^1_t\notempty\,\,\forall
    t\geq0\right\}.\] By \eqref{eq:bigM}, $\pp(\mathrm{G}_3)>1-\varepsilon$. Observe that
  on the event $\mathrm{G}_1\cap\mathrm{G}_2\cap\mathrm{G}_3$, $\eta^2_t$ contains 1's at
  all times with probability 1. Therefore
  \begin{align*}
    \pp^{\chi_{[-n,n]^d}}\!\left(A_t\notempty\,\,\forall
      t\geq0\right)&\geq\pp(\mathrm{G}_1\cap\mathrm{G}_2\cap\mathrm{G}_3)\\
    &\geq1-\pp(\mathrm{G}_1^\text{c})-\pp(\mathrm{G}_2^\text{c})-\pp(\mathrm{G}_3^\text{c})\\
    &>1-3\varepsilon,
  \end{align*}
  and choosing $\varepsilon$ small enough we get the result.
\end{proof}

In the following lemma we combine and extend for our process the results in \citet{ligg2}
that lead to the block conditions. Consider the process ${}_L\eta_t$, for $L>0$, where no
births are allowed outside of $(-L,L)^d$. Define $N_+(L,T)$ to be the maximal number of
space-time points in
\[S_+(L,T)=\{(x,s)\in(\{L\}\times[0,L)^{d-1})\times[0,T]\!:\,x\in{}_LA_s\}\] such that
each pair of these points having the same spatial coordinate have their time coordinates
at distance at least 1.

\begin{lem}\label{lem:techBlock}
  Suppose that the process survives. Then for any $\sigma>0$ there is a positive integer
  $n$ satisfying the following: for any given pair of positive integers $N$ and $M$, there
  are choices of a positive integer $L$ and a positive real number $T$ such that
  \begin{subequations}\label{eq:techBlock}
    \begin{gather}
      \pp^{\chi_{[-n,n]^d}}\!\left(\left|{}_LA_T\cap[0,L)^d\right|>N\right)\geq1-\sigma^{2^{-d}}\label{eq:techBlock1}\\
      \intertext{and}
      \pp^{\chi_{[-n,n]^d}}\!\left(\vphantom{2^2}N_+(L,T)>M\right)\geq1-\sigma^{2^{-d}/d}.\label{eq:techBlock2}
    \end{gather}
  \end{subequations}
\end{lem}

\begin{proof}
  By Lemma \ref{lem:bigN} we can choose a large enough integer $n$ such that
  \begin{equation}\label{eq:bigN}
    \pp^{\chi_{[-n,n]^d}}\!\left(A_t\notempty\,\,\forall
      t\geq0\right)>1-\sigma^2.
  \end{equation}
  Having this, the proof of the lemma is a simple adaptation of the corresponding proofs
  for the contact process. To avoid repetition of published results, we will explain the
  main ideas involved and why the original proofs still work with the random environment,
  but refer the reader to Section I.2 of \citet{ligg2} for the details.

  We claim the following: for any finite $A\subseteq\zz^d$ and any $N\geq1$,
  \begin{equation}\label{eq:prop2.2}
    \lim_{t\rightarrow\infty}\lim_{L\rightarrow\infty}\pp^{\chi_A}\!\left(\vphantom{2^2}\left|{}_LA_t\right|\geq
      N\right)=\pp^{\chi_A}\!\left(\vphantom{2^2}A_t\notempty\,\,\forall
      t\geq0\right).
  \end{equation}
  To see that this is true, we observe that
  \[\lim_{L\rightarrow\infty}\pp^{\chi_A}\!\left(\vphantom{2^2}\left|{}_LA_t\right|\geq
    N\right)=\pp^{\chi_A}\!\left(\vphantom{2^2}\left|A_t\right|\geq N\right)\] and then
  argue that, conditioned on survival, $|A_t|\longrightarrow\infty$ as
  $t\longrightarrow\infty$ with probability 1. This follows from the easy fact that there
  is an $\varepsilon_N>0$ such that if $|A|\leq N$ then the process started with 1's at
  $A$ becomes extinct with probability at least $\varepsilon_N$, so
  \[\pp^{\chi_A}\!\left(0<|A_t|\leq
    N\right)\varepsilon_N\leq\pp^{\chi_A}\!\left(t<\tau<\infty\right)\xrightarrow[t\rightarrow\infty]{}0.\]

  The next step is to use positive correlations to localize estimates on the cardinality
  of ${}_LA_t$ to a specific orthant of $\zz^d$: for every $N\geq1$ and $L\geq n$,
  \begin{equation}\label{eq:prop2.6}
    \pp^{\chi_{[-n,n]^d}}\!\left(\vphantom{2^2}\left|{}_LA_t\cap[0,L)^d\right|\leq
      N\right)\leq\left[\pp^{\chi_{[-n,n]^d}}\!\left(\vphantom{2^2}\left|{}_LA_t\right|\leq
        2^dN\right)\right]^{2^{-d}}.
  \end{equation}
  This relation follows easily from the positive correlations result in Lemma
  \ref{lem:posCorr}, and its proof is identical the proof of Proposition I.2.6 in
  \citet{ligg2}.

  Observe that \eqref{eq:bigN}, \eqref{eq:prop2.2}, and \eqref{eq:prop2.6} together
  suffice to obtain \eqref{eq:techBlock1}. The preceding arguments can be modified to
  obtain similar estimates for $N_+(L,T)$, which in turn give \eqref{eq:techBlock2}. The
  only detail remaining is getting the same $L$ and $T$ to work for both
  inequalities. This is done by obtaining sequences $L_j\nearrow\infty$ and
  $T_j\nearrow\infty$ such that \eqref{eq:techBlock1} holds with $L=L_j$ and $T=T_j$ for
  every $j\geq1$, and then adapting the arguments above to show that \eqref{eq:techBlock2}
  must hold for some pair $(L_j,T_j)$.  We refer the reader to the proof of Theorem I.2.12
  in Liggett's book for the details on how this is achieved, and remark that the argument
  depends only on properties such as positive correlations and the Feller property which
  are available both for $\eta_t$ and the contact process.
\end{proof}

We state now the block conditions that are equivalent to the survival of the process.

\begin{teosec}\label{thm:block}
  The process survives if and only if for any given $\varepsilon>0$ there are positive
  integers $n$ and $L$ and a positive real number $T$ such that the following conditions
  \hypertarget{BC}{{\rm(BC)}} are satisfied:
  \begin{multline}
    \label{eq:surv1}\tag{BC1}
    \pp^{\chi_{[-n,n]^d}}\!\left({}_{L+2n}A_{T+1}\supseteq x+[-n,n]^d\text{ for some }x\in[0,L)^d\right)\\
    >1-\varepsilon
  \end{multline}
  and
  \begin{multline}
    \label{eq:surv2}\tag{BC2}
    \pp^{\chi_{[-n,n]^d}}\!\left(\vphantom{2^2}{}_{L+2n}A_{t+1}\supseteq x+[-n,n]^d\text{ for some }0\leq t\leq T\right.\\
    \left.\vphantom{2^2}\text{ and some
      }x\in\{L+n\}\times[0,L)^{d-1}\right)>1-\varepsilon.
  \end{multline}
\end{teosec}

Observe that these conditions correspond exactly to the conditions in Theorem I.2.12 of
\citet{ligg2}. This will allow us to borrow the arguments from Liggett's book to prove
that {\BC} implies survival for $\eta_t$. The reason why we need the conditions {\BC}
starting $\eta_t$ from $\chi_{[-n,n]^d}$ is because the proof of their sufficiency for
survival (as well as their use in the proof of Theorem \ref{thm:compConv}) demands
obtaining repeatedly cubes fully occupied by 1's and, at each step, restarting the process
at the lowest possible configuration having those cubes fully occupied.

\begin{proof}[Proof of Theorem \ref{thm:block}]
  The proof uses the exact same arguments as those in the proofs of Theorems I.2.12 and
  I.2.23 in \citet{ligg2}. As before, we will only make some remarks and refer the reader
  to Liggett's book for the details.

  The necessity of {\BC} follows from Lemma \ref{lem:techBlock}, by choosing the
  quantities $N$ and $M$ to be large enough to produce the desired boxes filled with 1's.

  For the sufficiency of {\BC}, attractiveness and \eqref{S2} imply that it is enough to
  show that for some $n>0$ the process started at $\chi_{[-n,n]^d}$ contains 1's at all
  times (by using, as above, the fact that for any given $n>0$ the process started at
  $\chi_{\{0\}}$ has $[-n,n]^d$ fully occupied by time 1 with some positive probability).
  The proof of this fact relies on starting with a large enough cube fully occupied by 1's
  and then moving its center in an appropriate way. This is used to compare the process
  with supercritical oriented site percolation, and conclude that such boxes exist for all
  times with positive probability.
\end{proof}

The following consequence of Theorem \ref{thm:block} is obtained in the same way as for
the contact process, see Theorem I.2.25 in \citet{ligg2} for the details.

\begin{cor}\label{cor:critical}
  If $\beta=\beta_c(\alpha,\delta)$ or $\delta=\delta_c(\alpha,\beta)$, then the process
  dies out.
\end{cor}

\section{Complete convergence}\label{sec:compConv}

We are ready now to use the block construction of Section \ref{sec:block} to prove Theorem
\ref{thm:compConv}. The key step in the proof will be to obtain the result in the special
case where the initial distribution $\mu$ is a probability measure of the form $\nA$, in
which case we can use duality.

\begin{prop}\label{prop:compConvNuA}
  For every $A\subseteq\zz^d$,
  \[\eta^{\nA}_t\Longrightarrow\pp^{\nA}\!\left(\tau<\infty\right)\unu
  +\pp^{\nA}\!\left(\tau=\infty\right)\onu.\]
\end{prop}

To prove the proposition we need a preliminary lemma. Both the proof of the proposition
and this lemma are inspired by the proof Theorem 2 in \citet{durrettMoller}.

We will denote by $\pp^{\nA,\nC}$ the probability measure associated to starting the
process at $\nA$ and its dual at $\nC$, using the same realization of the graphical
representation, as explained in Section \ref{subsec:duality}.

\begin{lem}\label{lem:dec}
  For every finite $C\subseteq\zz^d$ and every $\varepsilon>0$, if $r$ is a positive real
  number and $s$ is large enough, then
  \begin{multline*}
    \left|\pp^{\nu_A,\nC}\!\left(\tau>\frac{s}{2},
        \dual{A}{r+s}_r\notempty,\dual{B}{r+s}_0\cap
        D\notempty\right)\right.\\
    \left.-\pp^{\nu_A}\!\left(\tau>\frac{s}{2}\right)
      \pp^\nC\!\left(\dual{A}{r}_r\notempty,\dual{B}{r}_0\cap
        D\notempty\right)\right|<\varepsilon.
  \end{multline*}
\end{lem}

Observe that for the (ordinary) contact process, the forward process and the dual are
independent when they run on nonoverlapping time intervals, so this fact is trivial and
holds with $s/2$ replaced by $s$.

\begin{proof}[Proof of Lemma \ref{lem:dec}]
  Given $r$ and $\varepsilon$, there is a $q=q(|C|)$ such that every dual active path in
  $(\dual{\eta}{\nC,r}_u)_{0\leq u\leq r}$ stays inside $C+[-q,q]^d$ with probability at
  least $1-\varepsilon$. To see this, observe that the number of particles in all such
  dual active paths is dominated by $X_r$, where $(X_r)_{r\geq0}$ is a branching process
  starting with $|C|$ particles and with birth rate $\beta$ and death rate $0$ (we are
  ignoring deaths and coalescence of paths).  By Markov's inequality,
  $\pp(X_r>q)\leq\ee(X_r)/q\leq\varepsilon$ for large enough $q$.  Since any dual active
  path in $\dual{\eta}{\nC,r}_t$ starts inside $C$, $X_r\leq q$ implies that all dual
  active paths are contained inside $C+[-q,q]^d$ up to time $r$.

  Now denote by $\eta^{(\mr,s/2)}_t$ and $\dual{\eta}{(\mr,s/2),r}_t$ modifications of the
  process and its dual, constructed on the same graphical representation as the original
  ones, where the environment is reset at time $s/2$ to its equilibrium $\mr$,
  independently of its state before $s/2$ (that is, at time $s/2$ we replace every $-1$ by
  a $0$ and then flip every site to $-1$ with probability $\rho$, regardless of it being
  at state 0 or 1). Then for given $r$ and $q$, if $s$ is large enough we have that
  \begin{multline}\label{eq:bdCpdEnvs}
    \pp^{\nA,\nC}\!\left(B_u=B^{(\mr,s/2)}_u\text{ on }C+[-q,q]^d\,\,\,\forall u\in[s,s+r]\right)\\
    \geq\left(1-e^{-\alpha(1+\delta)s/2}\right)^{\left|C+[-q,q]^d\right|} >1-\varepsilon.
  \end{multline}
  Indeed, for any given $x\in C+[-q,q]^d$ the probability that $B_u$ and $B^{(\mr,s/2)}_u$
  ar equal at $x$ for every $u\in[s,s+r]$ is bounded below by the probability that an
  exponential random variable with parameter $\alpha(1+\delta)$ is smaller than $s/2$
  (because any symbol $\bullet_{-1}$ or $\ast_{-1}$ above $(x,s/2)$ leaves the environment
  at that site equal for both processes from that time on).

  The property discussed at the first paragraph of the proof together with
  \eqref{eq:bdCpdEnvs} imply that
  \begin{multline*}
    \left|\pp^{\nA,\nC}\!\left(\tau>\frac{s}{2},
        \dual{A}{r+s}_r\notempty,\dual{B}{r+s}_0\cap
        D\notempty\right)\right.\\
    \left.-\pp^{\nA,\nC}\!\left(\tau>\frac{s}{2},\dual{A}{(\mr,s/2),r+s}_r\notempty,\dual{B}{(\mr,s/2),r+s}_0\cap
        D\notempty\right)\right|<2\varepsilon.
  \end{multline*}
  The statement of the lemma follows now from the independence of disjoint parts of the
  graphical representation and the stationarity of $B_t$, since
  \begin{align*}
    \pp^{\nA,\nC}\!\left(\vphantom{2^2}\tau>\frac{s}{2},\right.
    &\left.\vphantom{2^2}\dual{A}{(\mr,s/2),r+s}_r\notempty,\dual{B}{(\mr,s/2),r+s}_0\cap
      D\notempty\right)\\
    &=\pp^{\nA}\!\left(\tau>\frac{s}{2}\right)
    \pp^{\nC}\!\left(\dual{A}{(\mr,s/2),r+s}_r\notempty,\dual{B}{(\mr,s/2),r+s}_0\cap
      D\notempty\right)\\
    &=\pp^{\nA}\!\left(\tau>\frac{s}{2}\right)\pp^{\nC}\!\left(\dual{A}{r}_r\notempty,\dual{B}{r}_0\cap
      D\notempty\right).\qedhere
  \end{align*}
\end{proof}

\begin{proof}[Proof of Proposition \ref{prop:compConvNuA}]
  The result is straightforward in the subcritical case. If the process survives, and
  since weak convergence in this setting corresponds to the convergence of the
  finite-dimensional distributions, it is enough to prove that the following three
  properties hold for any two finite subsets $C,D$ of $\zz^d$:
  \begin{gather}
    \label{eq:c1} \pp^\nA\!\left(A_t\cap
      C\notempty\right)\xrightarrow[t\rightarrow\infty]{}\pp^\nA\!\left(\tau=\infty\right)\onu\left(\vphantom{2^2}
      \left\{(E,F)\!:\,E\cap
        C\notempty\right\}\right),\tag{c1}\\
    \label{eq:c2}
    \begin{split}
      \pp^\nA\!\left(B_t\cap
        D\notempty\right)&=\pp^\nA\left(\tau<\infty\right)\unu\left(\vphantom{2^2}\left\{(E,F)\!:\,F\cap
          D\notempty\right\}\right)\\
      &\hspace{0.65in}+\pp^\nA\!\left(\tau=\infty\right)\onu\left(\vphantom{2^2}\left\{(E,F)\!:\,F\cap
          D\notempty\right\}\right),
    \end{split}
    \tag{c2}\\
    \intertext{and}
    \label{eq:c3}
    \begin{split}
      &\pp^\nA\!\left(A_t\cap C\notempty,B_t\cap
        D\notempty\right)\\
      &\hspace{0.6in}\xrightarrow[t\rightarrow\infty]{}
      \pp^\nA\!\left(\tau=\infty\right)\onu\left(\vphantom{2^2}\left\{(E,F)\!:\,E\cap
          C\notempty,F\cap D\notempty\right\}\right).
    \end{split}
    \tag{c3}
  \end{gather}
  Indeed, all the finite-dimensional distributions of the process are determined by these
  probabilities via the inclusion-exclusion formula. Observe that the right side of
  \eqref{eq:c2} is equal to $\mr\!\left(\left\{\eta\!:\,\eta(x)=-1 \text{ for some }x\in
      D\right\}\right)$.

  \eqref{eq:c1} follows from the same arguments used in \cite{ligg2} for the contact
  process. Using duality (Proposition \ref{prop:dual}), the proof of Theorem I.1.12 in
  that book applies in the same way to obtain the fact that \eqref{eq:c1} holds if and
  only if for every $x\in\zz^d$ and every $A\subseteq\zz^d$,
  \begin{subequations}
    \begin{gather}
      \pp^\nA\!\left(\tau=\infty\right)=\pp^\nA\!\left(x\in A_t\,\text{ i.o.}\right)\label{eq:ccl1}\\
      \intertext{and}
      \lim_{n\rightarrow\infty}\liminf_{t\rightarrow\infty}\pp^{\nu_{[-n,n]^d}}\!\left(A_t\cap[-n,n]^d\notempty\right)=1.\label{eq:ccl2}
    \end{gather}
  \end{subequations}
  The analogous conditions are checked for the contact process in the proof of Theorem
  I.2.27 in \cite{ligg2}. \eqref{eq:ccl1} follows from the same proof after some minor
  modifications, so we will skip the argument.  For \eqref{eq:ccl2}, Theorem
  \ref{thm:block} allows us to use Liggett's arguments to get the desired limit when
  $\nu_{[-n,n]^d}$ is replaced by $\chi_{[-n,n]^d}$, so given any $\varepsilon>0$ we can
  choose a large enough integer $m$ such that
  \begin{equation}\label{eq:bigMCcl2}
    \liminf_{t\rightarrow\infty}\pp^{\chi_{[-m,m]^d}}\!\left(A_t\cap[-m,m]^d\notempty\right)>1-\varepsilon.
  \end{equation}
  Given this $m$, we can choose a large enough $n$ so that the process started at
  $\nu_{[-n,n]^d}$ contains at time 0 a fully occupied cube of side $2m+1$ (contained in
  $[-n,n]^d$) with probability at least $1-\varepsilon$ (in fact, any translate of
  $[-m,m]^d$ contained in $[-n,n]^d$ is fully occupied by 1's with some probability $p>0$,
  so we only need to choose $n$ so that $[-n,n]^d$ contains enough disjoint translates of
  $[-m,m]^d$). 
  On this event we can restart the
  process by putting every site outside that cube at state $-1$ and use attractiveness,
  translation invariance, and \eqref{eq:bigMCcl2} to get
  \[\liminf_{t\rightarrow\infty}\pp^{\nu_{[-n,n]^d}}\!\left(A_t\cap[-n,n]^d\notempty\right)>(1-\varepsilon)^2,\]
  whence \eqref{eq:ccl2} follows. There is only one detail to consider: in his book,
  Liggett only proves the condition analogous to \eqref{eq:ccl2} in the case $d\geq2$,
  because it is simpler and the case $d=1$ was already done in \citet{ligg1}. The
  difficulty in the one-dimensional case arises from the fact that certain block events
  are not independent.  This can be overcome by comparing with $k$-dependent oriented site
  percolation instead of ordinary oriented site percolation (see Theorem B26 in
  \cite{ligg2}). We refer the reader to Section 5 of \cite{durretSchonmann}, where the
  authors use a similar block construction to derive a complete convergence theorem for a
  general class of one-dimensional growth models.

  \eqref{eq:c2} is trivial due to the stationarity of the environment process. To prove
  \eqref{eq:c3} we start by observing that
  \begin{multline}\label{eq:rwProb}
    \pp^\nA\!\left(A_{r+s}\cap C\notempty,B_{r+s}\cap D\notempty\right)\\
    =\pp^{\nA,\nC}\!\left(A_s\cap\dual{A}{r+s}_r\notempty,\dual{B}{r+s}_0\cap
      D\notempty\right),
  \end{multline}
  which follows from constructing $(\eta^\nA_u)_{0\leq u\leq r+s}$ and
  $(\dual{\eta}{\nC,r+s}_u)_{0\leq u\leq r+s}$ on the same copy of the graphical
  representation. On the other hand,
  \begin{equation}\label{eq:noB}
    \begin{aligned}
      &\left|\pp^{\nA,\nC}\!\left(A_s\cap\dual{A}{r+s}_r\notempty,\dual{B}{r+s}_0\cap D\notempty\right)\right.\\
      &\hspace{0.55in}\left.\hphantom{\pp^{\nA,\nC} A_s}-\pp^{\nA,\nC}\!\left(A_s\notempty,\dual{A}{r+s}_r\notempty,\dual{B}{r+s}_0\cap D\notempty\right)\right|\\
      &\hspace{0.4in}=\pp^{\nA,\nC}\!\left(A_s\notempty,\dual{A}{r+s}_r\notempty,A_s\cap\dual{A}{r+s}_r=\emptyset,\dual{B}{r+s}_0\cap D\notempty\right)\\
      &\hspace{0.4in}\leq\pp^{\nA,\nC}\!\left(A_s\notempty,\dual{A}{r+s}_r\notempty,A_s\cap\dual{A}{r+s}_r=\emptyset\right)\\
      &\hspace{0.4in}=\pp^{\nA,\nC}\!\left(A_s\notempty,\dual{A}{r+s}_r\notempty\right)-\pp^{\nA,\nC}\!\left(A_s\cap\dual{A}{
          r+s}_r\notempty\right).
    \end{aligned}
  \end{equation}

  Observe that
  \begin{equation*}
    \pp^\nA\!\left(s/2<\tau<\infty\right)\xrightarrow[s\rightarrow\infty]{}0.
  \end{equation*}
  Thus, for any given $D\subseteq\zz^d$ and $\varepsilon>0$, and for large enough $s$, we
  can write
  \begin{equation}\label{eq:largeTau}
    \begin{aligned}
      &\left|\pp^{\nA,\nC}\!\left(A_s\notempty,\dual{A}{r+s}_r\notempty,\dual{B}{r+s}_0\cap D\notempty\right)\right.\\
      &\hspace{1.5in}\left.-\pp^{\nA,\nC}\!\left(\tau>s/2,\dual{A}{r+s}_r\notempty,\dual{B}{r+s}_0\cap D\notempty\right)\right|\\
      &\hspace{0.5in}=\pp^{\nA,\nC}\!\left(s/2<\tau\leq s,\dual{A}{r+s}_r\notempty,\dual{B}{r+s}_0\cap D\notempty\right)\\
      &\hspace{0.5in}\leq\pp^\nA\!\left(s/2<\tau<\infty\right)<\frac{\varepsilon}{3}.
    \end{aligned}
  \end{equation}
  
  Putting the previous observations together we get, for large enough $s$
  \begin{align*}
    &\left|\pp^\nA\!\left(\vphantom{\dual{A}{r+s}}A_{r+s}\cap C\notempty,B_{r+s}\cap D\notempty\right)\right.&&\\
    &\hspace{0.8in}\left.-\pp^{\nA,\nC}\!\left(A_s\notempty,\dual{A}{r+s}_r\notempty,\dual{B}{r+s}_0\cap D\notempty\right)\right|&&\\
    &\hspace{0.12in}=\left|\pp^{\nA,\nC}\!\left(A_s\cap\dual{A}{r+s}_r\notempty,\dual{B}{r+s}_0\cap D\notempty\right)\right.&&\\
    &\hspace{0.8in}\left.-\pp^{\nA,\nC}\!\left(A_s\notempty,\dual{A}{r+s}_r\notempty,\dual{B}{r+s}_0\cap D\notempty\right)\right|&&\text{by \eqref{eq:rwProb}}\\
    &\hspace{0.12in}\leq\pp^{\nA,\nC}\!\left(A_s\notempty,\dual{A}{r+s}_r\notempty\right)-\pp^{\nA,\nC}\!\left(A_s\cap\dual{A}{r+s}_r\notempty\right)&&\text{by \eqref{eq:noB}}\\
    &\hspace{0.12in}\leq\frac{\varepsilon}{3}+\left|\pp^{\nA,\nC}\!\left(\tau>s/2,\dual{A}{r+s}_r\notempty\right)-\pp^{\nA,\nC}\!\left(A_s\cap\dual{A}{r+s}_r\notempty\right)\right|&&\text{by
      \eqref{eq:largeTau}},
  \end{align*}
  where we used $D=\zz^d$ and the fact that $\pp^{\nA,\nC}(\dual{B}{r+s}_0\notempty)=1$ in
  the application of \eqref{eq:largeTau}. Using again this fact to apply Lemma
  \ref{lem:dec} with $D=\zz^d$, and then using duality we get
  \begin{align*}
    &\left|\pp^{\nA,\nC}\!\left(\tau>s/2,\dual{A}{r+s}_r\notempty\right)-\pp^{\nA,\nC}\!\left(A_s\cap\dual{A}{r+s}_r\notempty\right)\right|\\
    &\hspace{0.2in}\leq\frac{\varepsilon}{3}+\left|\pp^\nA\!\left(\tau>s/2\right)\pp^\nC\!\left(\vphantom{\dual{A}{r+s}}\dual{A}{r}_r\notempty\right)-\pp^{\nA,\nC}\!\left(A_s\cap\dual{A}{r+s}_r\notempty\right)\right|\\
    &\hspace{0.2in}=\frac{\varepsilon}{3}+\left|\vphantom{\dual{A}{r+s}}\pp^\nA\!\left(\tau>s/2\right)\pp^{\nu_{\zz^d}}\!\left(A_r\cap
        C\notempty\right)-\pp^\nA\!\left(A_{r+s}\cap C\notempty\right)\right|
  \end{align*}
  for large enough $s$. By \eqref{eq:c1}, the last difference converges to 0 as
  $r,s\rightarrow\infty$, so we finally get
  \begin{multline*}
    \left|\pp^\nA\!\left(\vphantom{\dual{A}{r+s}}A_{r+s}\cap C\notempty,B_{r+s}\cap D\notempty\right)\right.\\
    \left.-\pp^{\nA,\nC}\!\left(A_s\notempty,\dual{A}{r+s}_r\notempty,\dual{B}{r+s}_0\cap
        D\notempty\right)\right|<\varepsilon
  \end{multline*}
  for large enough $r,s$.

  This calculation implies that in order to prove \eqref{eq:c3} it is enough to show that
  \begin{multline*}
    \pp^{\nA,\nC}\!\left(A_s\notempty,\dual{A}{r+s}_r\notempty,\dual{B}{r+s}_0\cap
      D\notempty\right)\\
    \xrightarrow[r,s\rightarrow\infty]{}\pp^\nA\!\left(\tau=\infty\right)
    \onu\left(\vphantom{2^2}(E,F)\!:\,E\cap C\notempty,F\cap D\notempty\right).
  \end{multline*}
  Repeating the previous application of \eqref{eq:largeTau} and Lemma \ref{lem:dec} we get
  that, for large enough $s$,
  \begin{align*}
    &\left|\pp^{\nA,\nC}\!\left(A_s\notempty,\dual{A}{r+s}_r\notempty,\dual{B}{r+s}_0\cap
        D\notempty\right)\right.\\
    &\hspace{1.2in}\left.-\pp^\nA\!\left(\tau>s/2\right)\pp^\nC\!\left(\dual{A}{r}_r\notempty,\dual{B}{r}_0\cap
        D\notempty\right)\right|\\
    &\hspace{0.2in}\leq\frac{\varepsilon}{2}+\left|\pp^{\nA,\nC}\!\left(\tau>s/2,\dual{A}{r+s}_r\notempty,\dual{B}{r+s}_0\cap
        D\notempty\right)\right.\\
    &\hspace{1.2in}\left.-\pp^\nA\!\left(\tau>s/2\right)\pp^\nC\!\left(\dual{A}{r}_r\notempty,\dual{B}{r}_0\cap
        D\notempty\right)\right|\\
    &\hspace{0.2in}\leq\varepsilon.
  \end{align*}
  Therefore, we can finally reduce to proving that
  \begin{multline*}
    \pp^\nA\!\left(\tau>s/2\right)\pp^\nC\!\left(\dual{A}{r}_r\notempty,\dual{B}{r}_0\cap
      D\notempty\right)\\
    \xrightarrow[r,s\rightarrow\infty]{}
    \pp^\nA\!\left(\tau=\infty\right)\onu\left(\vphantom{2^2}(E,F)\!:\,E\cap
      C\notempty,F\cap D\notempty\right).
  \end{multline*}
  This follows easily from duality, since \eqref{eq:dual1} yields
  \begin{multline*}
    \pp^\nA\!\left(\tau>s/2\right)\pp^\nC\!\left(\dual{A}{r}_r\notempty,\dual{B}{r}_0\cap D\notempty\right)\\
    =\pp^\nA\!\left(\tau>s/2\right)\pp^{\nu_{\zz^d}}\!\left(\vphantom{2^2}A_r\cap
      C\notempty,B_r\cap D\notempty\right),
  \end{multline*}
  and this last term converges to the desired limit as $r,s\rightarrow\infty$.
\end{proof}

We extend now Proposition \ref{prop:compConvNuA} to the general case.

\begin{proof}[Proof of Theorem \ref{thm:compConv}]
  It is enough to show that
  \begin{equation}\label{eq:weakConv}
    \lim_{t\rightarrow\infty}\ee^\mu(f(\eta_t))
    =\pp^{\mu}(\tau<\infty)\int\!f\,d\unu+\pp^{\mu}(\tau=\infty)\int\!f\,d\onu
  \end{equation}
  for every $f$ in the space of continuous increasing functions depending on finitely many
  coordinates of $\cx$, which we will denote by $\cf$. To see this, observe that given any
  two finite subsets $C,D$ of $\zz^d$, the functions
  \[f_1(E,F)=\uno{E\cap C\notempty},\,\,f_2(E,F)=\uno{F\cap D=\emptyset},\,\,\text{and }
  f_3(E,F)=\uno{E\cap C\notempty,F\cap D=\emptyset}\] are all in $\cf$ and (as in the proof of
  Proposition \ref{prop:compConvNuA}) all the finite-dimensional distributions of the
  process can be obtained from $\ee^\mu(f_1(\eta_t))$, $\ee^\mu(f_2(\eta_t))$, and
  $\ee^\mu(f_3(\eta_t))$ by the inclusion-exclusion formula.
  
  Let $f$ be a function in $\cf$ and observe that, in particular, $f$
  is bounded. One inequality in \eqref{eq:weakConv} is easy: by the Markov property and
  attractiveness, given $0<s<t$ we have that
  \begin{align*}
    \ee^\mu(f(\eta_t))&=\ee^\mu(f(\eta_t),\,\tau<s)
    +\ee^\mu(f(\eta_t),\,\tau\geq s)\\
    &=\ee^\mu(\ee^{\eta_s}(f(\eta_{t-s}),\,\tau<s)
    +\ee^\mu(\ee^{\eta_s}(f(\eta_{t-s})),\,\tau\geq s)\\
    &\leq\ee(f(\eta^0_{t-s}))\pp^\mu(\tau<s)
    +\ee^{\chi_{\zz^d}}(f(\eta_{t-s}))\pp^\mu(\tau\geq s),
  \end{align*}
  where $\eta^0_t$ denotes the process started at the configuration $\eta\equiv0$. Since
  $\eta^0_t\Longrightarrow\mr=\unu$ and $\eta^{\chi_{\zz^d}}_t\Longrightarrow\onu$, we get
  \[\limsup_{t\rightarrow\infty}\ee^\mu(f(\eta_t))
  \leq\pp^\mu(\tau<s)\int\!f\,d\unu+\pp^\mu(\tau\geq s)\int\!f\,d\onu,\] and now taking
  $s\rightarrow\infty$ we deduce that
  \begin{equation}\label{eq:strMarkovSup}
    \limsup_{t\rightarrow\infty}\ee^\mu(f(\eta_t))
  \leq\pp^\mu(\tau<\infty)\int\!f\,d\unu+\pp^\mu(\tau=\infty)\int\!f\,d\onu.
\end{equation}

  To obtain the other inequality in \eqref{eq:weakConv} we will begin by considering the
  case $\mu=\chi_{[-n,n]^d}$ and showing that, given any $\varepsilon>0$ and any $x\in\zz^d$,
  \begin{equation}\label{eq:bdIntChi}
    \liminf_{t\rightarrow\infty}\ee^{\chi_{x+[-n,n]^d}}(f(\eta_t),\,\tau=\infty)
    \geq\int\!f\,d\onu-\varepsilon
  \end{equation}
  for large enough $n$. By the translation invariance of $\eta_t$ and $\onu$, it is enough
  to consider the case $x=0$. To show \eqref{eq:bdIntChi} we will use the construction
  introduced in the proof of Lemma \ref{lem:bigN}. Using the notation of that proof,
  recall that we showed that, given any $\gamma>0$, there are positive integers $n>k>m$
  such that
  \[\pp(\mathrm{G}_1\cap\mathrm{G}_2\cap\mathrm{G}_3)>1-3\gamma.\]
  This means that the processes $\eta^1_t$ (started at $\nu_{[-m,m]^d}$) and $\eta^2_t$
  (started at $\chi_{[-n,n]^d}$) can be coupled in such a way that, with probability at
  least $1-3\gamma$, for all $t\geq0$ we have that $A^1_t\notempty$,
  $A^2_t\notempty$, $A^1_t\subseteq Q(k+1+ct)$, and $B^2_t\subseteq B^1_t$ inside $Q(k+1+ct)$.

  Let $\mathrm{G}=\mathrm{G}_1\cap\mathrm{G}_2\cap\mathrm{G}_3$ and $\gamma>0$ and choose
  $n>k>m$ so that $\pp(\mathrm{G})>1-3\gamma$. We will denote by $\tau^1$ and
  $\tau^2$ the extinction times of the processes $\eta^1_t$ and $\eta^2_t$, respectively.
  Define
  \[\check{\eta}_t=(A^1_t,B^2_t)\]
  and observe that, on the event $\mathrm{G}$, $\check{\eta}_t$ defines an $\cx$-valued
  process and, moreover, $\eta^2_t\geq\check{\eta}_t$ for all $t\geq0$. Therefore, since
  $f$ is increasing and $\{\tau^2=\infty\}\supseteq\mathrm{G}$,
  \begin{equation}\label{eq:limEta2}
    \ee(f(\eta^2_t),\,\tau^2=\infty)
    \geq\ee(f(\eta^2_t),\mathrm{G})
    \geq\ee(f(\check{\eta}_t),\mathrm{G})
  \end{equation}
  for all $t\geq0$. Now observe that, trivially,
  \[\ee(f(\check{\eta}_t),\mathrm{G})=\ee(f(\check{\eta}_t),\,\tau^2=\infty)
  -\ee(f(\check{\eta}_t),\tau^2=\infty,\mathrm{G}^\text{c}),\] and
  \[\ee(f(\check{\eta}_t),\tau^2=\infty,\mathrm{G}^\text{c})
  \leq\|f\|_\infty\pp(\mathrm{G}^\text{c})<3\gamma\|f\|_\infty,\] so
  \begin{equation}\label{eq:bdEtaCh1}
    \ee(f(\check{\eta}_t),\mathrm{G})
    >\ee(f(\check{\eta}_t),\,\tau^2=\infty)-3\gamma\|f\|_\infty
  \end{equation}
  for all $t\geq0$. On the other hand,
  \begin{equation}\label{eq:bdEtaCh2}
    \left|\ee\big(f(\check{\eta}_t),\,\tau^2=\infty\big)
      -\ee\big(f(\eta^1_t),\tau^2=\infty\big)\right|
    \xrightarrow[t\rightarrow\infty]{}0.
  \end{equation}
  To see this, observe that since $f$ depends on finitely many coordinates, then given any
  $q>0$, $f(\check{\eta}_s)=f(\eta^1_s)$ for all $s\geq t$ with probability at least $1-q$
  if $t$ is large enough. Indeed, if $K\subseteq\zz^d$ is the finite set of coordinates of
  $\cx$ on which $f$ depends, then repeating the calculations that led to \eqref{eq:cone2}
  we get that
  \begin{multline*}
    \hspace{-0.1in}\pp\Big(B^1_s(x)\neq B^2_s(x)\text{ for some $x\in K$ and some $s\geq t$}\Big)\\
    \leq\sum_{x\in K}\pp\!\left(\text{no $\bullet_{-1}$ or $\ast_{-1}$ at $x$ by time
        $t$}\right)
    =|K|e^{-\alpha(1+\delta)t}\xrightarrow[t\rightarrow\infty]{}0.\hspace{-0.1in}
  \end{multline*}
  Therefore, given any $q>0$,
  \begin{equation*}
    \left|\ee\big(f(\check{\eta}_t\big),\,\tau^2=\infty)
      -\ee\big(f(\eta^1_t),\,\tau^2=\infty\big)\right|
    \leq\ee\!\left(\left|f(\check{\eta}_t)-f(\eta^1_t)\right|\right)
    \leq2q\|f\|_\infty
  \end{equation*}
  for large enough $t$, and we get \eqref{eq:bdEtaCh2}.  Finally, we have that
  \begin{align*}
    \ee(f(\eta^1_t)&,\,\tau^2=\infty)\\
    &=\ee(f(\eta^1_t),\,\tau^1=\infty) -\left(\ee(f(\eta^1_t),\,\tau^1=\infty)
      -\ee(f(\eta^1_t),\mathrm{G})\right)\\
    &\hspace{0.4in}-\left(\ee(f(\eta^1_t),\mathrm{G})
      -\ee(f(\eta^1_t),\,\tau^2=\infty\right),
  \end{align*}
  and since $\mathrm{G}\subseteq\{\tau^1=\infty\}\cap\{\tau^2=\infty\}$,
  \[\left|\ee(f(\eta^1_t),\,\tau^i=\infty)-\ee(f(\eta^1_t),\mathrm{G})\right|
  \leq\|f\|_\infty\pp(\mathrm{G}^\text{c})<3\gamma\|f\|_\infty\] for $i=1,2$. Thus
  Proposition \ref{prop:compConvNuA} implies that
  \[\liminf_{t\rightarrow\infty}\ee(f(\eta^1_t),\,\tau^2=\infty)
    >\pp(\tau^1=\infty)\int\!f\,d\onu-6\gamma\|f\|_\infty,\]
  and since $\pp(\tau^1=\infty)\geq\pp(\mathrm{G})>1-3\gamma$ we obtain
    \begin{equation}\label{eq:bdEtaCh3}
    \liminf_{t\rightarrow\infty}\ee(f(\eta^1_t),\,\tau^2=\infty)
    >\int\!f\,d\onu-9\gamma\|f\|_\infty.
  \end{equation}
  Putting \eqref{eq:limEta2}, \eqref{eq:bdEtaCh1}, \eqref{eq:bdEtaCh2}, and
   \eqref{eq:bdEtaCh3} together we deduce that
   \[\liminf_{t\rightarrow\infty}\ee(f(\eta^2_t),\,\tau^2=\infty)
   \geq\int\!f\,d\onu-12\gamma\|f\|_\infty,\]
  and choosing $\gamma$ appropriately we obtain
  \eqref{eq:bdIntChi}.

  Getting back to the proof of the remaining inequality in \eqref{eq:weakConv}, let
  $\varepsilon>0$ and choose $n\in\nn$ so that \eqref{eq:bdIntChi} holds. Define
  \[N=\inf\left\{k\in\nn\!:\eta_k\supseteq x+[-n,n]^d\,\text{ for some
      $x\in\zz^d$}\right\}\]
  and let $p=\pp^{\chi_{\{0\}}}\!\big(A_1\supseteq x+[-n,n]^d\text{ for some
  }x\in\zz^d\big)>0$.
  Observe that for any $k\geq0$, if $A_k\notempty$ then $A_{k+1}$ contains
  some translate of $[-n,n]^d$ with probability at least $p$ (by attractiveness and
  translation invariance) and therefore, since the Poisson processes used in the graphical
  representation for disjoint time intervals are independent, we deduce that
  \begin{equation}\label{eq:NTau}
  \{\tau=\infty\}\subseteq\{N<\infty\}.
  \end{equation}
  When $N<\infty$ we will denote by $X$ the center of the
  corresponding fully occupied box. If there is more than one point $x$ such that
  $x+[-n,n]^d$ is fully occupied by 1's at time $N$, we pick $X$ to be the one minimizing
  $\phi(x)$, where $\phi$ is any fixed bijection between $\zz^d$ and $\nn$ (this ensures
  that the events $\{X=x\}$ are disjoint for different $x$). Then given $m\in\nn$, the
  Markov property and attractiveness imply that
  \begin{align*}
    \ee^\mu(f(\eta_t),\,&\tau=\infty)
    \geq\sum_{k=0}^m\ee^\mu(f(\eta_t),\,\tau=\infty,\,N=k)\\
    &=\sum_{k=0}^m\ee^\mu\big(\ee^{\eta_{k}}(f(\eta_{t-k}),\,\tau=\infty),\,N=k\big)\\
    &\geq\sum_{k=0}^m\sum_{x\in\zz^d}\ee^\mu\!\left(\ee^{\chi_{x+[-n,n]^d}}(f(\eta_{t-k}\right),\,
    \tau=\infty),\,N=k\,,X=x)
  \end{align*}
  for $t\geq m$. Since $f$ is bounded, \eqref{eq:bdIntChi} implies that
  \begin{align*}
    \liminf_{t\rightarrow\infty}\ee^\mu(f(\eta_t),\,\tau=\infty)
    &\geq\left(\int\!f\,d\onu-\varepsilon\right)
    \sum_{k=0}^m\sum_{x\in\zz^d}\pp^\mu(N=k,\,X=x)\\
    &=\left(\int\!f\,d\onu-\varepsilon\right)\pp^\mu(N\leq m).
  \end{align*}
  Taking now $m\rightarrow\infty$ we get by \eqref{eq:NTau} that
  \begin{align*}
    \liminf_{t\rightarrow\infty}\ee^\mu(f(\eta_t),\,\tau=\infty)
  &\geq\left(\int\!f\,d\onu-\varepsilon\right)\pp^\mu(N<\infty)\\
  &\geq\left(\int\!f\,d\onu-\varepsilon\right)\pp^\mu(\tau=\infty)
\end{align*}
 if $\varepsilon<\int\!f\,d\onu$, and taking $\varepsilon\rightarrow0$ we deduce that
  \begin{equation*}
    \liminf_{t\rightarrow\infty}\ee^\mu(f(\eta_t),\,\tau=\infty)
    \geq\pp^\mu(\tau=\infty)\int\!f\,d\onu.
  \end{equation*}
  On the other hand, by arguments similar to those that led to \eqref{eq:strMarkovSup}
  (using attractiveness to compare with the process started at $\chi_\emptyset$) we get
  \[\liminf_{t\rightarrow\infty}\ee^\mu(f(\eta_t),\,\tau<\infty)\geq\pp^\mu(\tau<\infty)\int\!f\,d\unu.\]
  We finally deduce that
  \[\liminf_{t\rightarrow\infty}\ee^\mu(f(\eta_t))
  \geq\pp^\mu(\tau<\infty)\int\!f\,d\onu+\pp^\mu(\tau=\infty)\int\!f\,d\onu,\] and the
  proof is ready.
\end{proof}

\section*{Acknowledgments}

I am grateful to my advisor, Rick Durrett, for his guidance during the work in this paper.
I thank him especially for his careful reading of several versions of the manuscript and
the thorough comments and suggestions he provided, which led to a substantial improvement
of the paper. I also wish to thank Ted Cox for reviewing part of an earlier version of the
manuscript. Finally, I thank an anonymous referee for many detailed comments and
suggestions which greatly improved the exposition of the paper.\\[-4pt]

\noindent I am also grateful to Dong Yao for pointing out a mistake in the published version of this article.

\bibliographystyle{natbib} \bibliography{biblio}

\end{document}